[1994/12/01]

\documentclass[twoside,reqno,10pt]{amsart}
\usepackage{amsmath,amssymb,amsthm}

\numberwithin{equation}{section}

\newtheorem{Cor}{Corollary}
\newtheorem{Thm}{Theorem}
\newtheorem{Prop}{Proposition}
\newtheorem{Lem}{Lemma}

\newtheorem*{Conj}{Conjecture}

\theoremstyle{definition}
\newtheorem{Def}{Definition}

\theoremstyle{remark}
\newtheorem{Rem}{Remark}

\renewcommand\ge{\geqslant}
\renewcommand\le{\leqslant}
\let\tildeaccent=\~
\renewcommand\~[1]{\widetilde{#1}}

\def\<{\left<}
\def\>{\right>}
\def\({\ifmmode\left(\else\textup{(}\fi}
\def\){\ifmmode\right)\else\textup{)}\fi}

\def\const{\operatorname{const}}

\def\id{{\mathbf 1}}

\renewcommand\:{\colon}
\def\R{{\mathbb R}}
\def\C{{\mathbb C}}

\def\f{\varphi}
\let\ssm=\smallsetminus
\def\pd#1#2{\frac{\partial#1}{\partial#2}}

\def\L{{\boldsymbol\varLambda}}
\def\ix{\operatorname{i}_{_X}}
\def\l{\lambda}
\def\G{\varGamma}

\let\paragraph=\S
\def\secref#1{\paragraph\ref{#1}}
\def\S{\varSigma}

\begin{document}

\title[Bounded Brieskorn decomposition]
 { Bounded decomposition in the Brieskorn lattice and Pfaffian Picard--Fuchs
 systems for Abelian integrals }

\author
 {Sergei Yakovenko}

\thanks{The research was supported by the Israeli Science Foundation grant
no.~18-00/1.}

\address{Department of  Mathematics\\
 Weizmann Institute of Science\\
 P.O.B.~26, Rehovot 76100\\
 Israel}
\email{{\tt yakov@wisdom.weizmann.ac.il}
\endgraf{\it WWW page\/}:
{\tt http://www.wisdom.weizmann.ac.il/\char'176 yakov/index.html}}
\date{January 2002}


\begin{abstract}
We suggest an algorithm for derivation of the Picard--Fuchs system
of Pfaffian equations for Abelian integrals corresponding to
semiquasihomogeneous Hamiltonians. It is based on an effective
decomposition of polynomial forms in the Brieskorn lattice. The
construction allows for an explicit upper bound on the norms of
the polynomial coefficients, an important ingredient in studying
zeros of these integrals.
\end{abstract}

\maketitle

\section{Introduction}

Given a polynomial in two variables $f\in\R[x,y]$ and a polynomial
1-form $\omega$ on $\R^2$, how many isolated ovals $\delta$ on the
level curves $f=\const$ may satisfy the condition
$\oint_\delta\omega=0$? This is the long-standing
\emph{infinitesimal Hilbert problem}, see
\cite{arnold:quel-prob-sys-dyn}. The answer is to be given in
terms of the degrees of $f$ and $\omega$.

A recent approach to this problem, suggested in
\cite{era-99,redundant,montreal} is based on the fact that periods
of polynomial 1-forms restricted on level curves of polynomials,
satisfy a system of differential equations with rational
coefficients, called the \emph{Picard--Fuchs system}. Under
certain restrictions on the monodromy group, the number of zeros
of solutions of such systems can be estimated from above in terms
of the magnitude of coefficients of this system, more precisely,
the norms of its matrix residues. Thus it becomes important to
derive the Picard--Fuchs system for Abelian integrals so
explicitly as to allow for the required estimates for the
residues.

In \cite{redundant} a Fuchsian system was derived in the
\emph{hypergeometric} form
\begin{equation}\label{hyperg}
  (t\cdot\id+A)\dot I=BI,\qquad \dot I=\tfrac{d}{dt}I(t),
\end{equation}
where $I(t)=(I_1(t),\dots,I_l(t))$ is a collection of integrals of
some monomial forms over any oval of the level curve $\{f=t\}$,
and $A,B$ are two constant $(l\times l)$-matrices of explicitly
bounded norms, depending on $f$ ($\id$ always stands for the
identity matrix of the appropriate size). The rational matrix
function $R(t)=(t\cdot\id+A)^{-1}B$ has only simple poles and the
norm of its matrix residues can be explicitly majorized provided
that the eigenvalues of $A$ remain well apart. This allows to
solve the infinitesimal Hilbert problem for all polynomials $f$
whose critical values (after a suitable normalization) are
sufficiently distant from each other. What remains is to study the
case of confluent critical values (including those at infinity).

In a general hypergeometric system \eqref{hyperg}, the residues
may or may not blow up as some of the singular points tend to each
other. The particular feature of the Picard--Fuchs system is its
\emph{isomonodromy}: the monodromy group remains the same under
deformations of $f$ (at least for sufficiently generic $f$). This
implies that even if the explosion of residues occurs, it cannot
be caused by the explosion of the eigenvalues. In order to find
out what indeed happens with the residues, the first step is to
write down as explicitly as possible the Picard--Fuchs system as a
flat meromorphic connexion with singularities in the holomorphic
bundle over the variety of all polynomials $f$ of a given degree.

This problem is solved in the paper for polynomials with a fixed
principal (quasi)homogeneous part having an isolated critical
point at the origin.

As an auxiliary first step, we need to describe explicitly the
structure of the relative cohomology module. While the subject is
fairly classic and sufficiently well understood, the existing
proofs do not allow for the quantitative analysis. We suggest an
alternative, completely elementary construction that immediately
yields all necessary bounds. This construction, exposed in
\secref{sec:relcom} is based on ``division by $f$'', a lemma
distilled from the paper \cite{francoise:brieskorn-thm} by
J.-P.~Fran\c{c}oise. The Pfaffian form of the Picard--Fuchs system
is derived in \secref{sec:picfuc}. In the last section we mention
some simple properties of the derived system and formulate a
conjecture that it has only logarithmic singularities in the
affine part.

\section{Relative cohomology revisited}\label{sec:relcom}

\subsection{Relative cohomology, Brieskorn and Petrov modules}
Denote by $\L^k$, $k=0,1,\dots,n$, the module of polynomial
$k$-forms on the complex affine space $\C^n$ for a fixed $n\ge 1$.
If $f\in\C[x_1,\dots,x_n]\simeq\L^0$ is a polynomial, then the
collection $df\land \L^{k-1}$ of $k$-forms \emph{divisible} by
$df\in\L^1$, is a $\C$-linear subspace in $\L^k$, and the quotient
\begin{equation}\label{relform}
    \L_f^k=\L^k/df\land \L^{k-1},\qquad k=1,\dots,n,
\end{equation}
is called the space of \emph{relative $k$-forms}. Since the
exterior derivative $d$ preserves divisibility by $df$, the
\emph{relative de~Rham complex} $\L^\bullet_f$,
\begin{equation}\label{rdR}
    0\longrightarrow\L^{1}_f\overset{d}\longrightarrow\L^{2}_f\cdots
    \overset{d}\longrightarrow\L^{n-1}_f\overset{d}\longrightarrow\L^n_f
    \overset{d}\longrightarrow0,
\end{equation}
naturally appears. A form $\omega\in\L^k$ is called
\emph{relatively closed} if $d\omega=df\land \eta$ and
\emph{relatively exact} if $\omega=df\land \xi+d\theta$ for
appropriate $\eta\in\L^k$ and $\xi,\theta\in\L^{k-1}$. The
\emph{relative cohomology groups} $\boldsymbol
H^k_f=H^k(\L^\bullet_f)$, relatively closed $k$-forms modulo
relatively exact ones, are important characteristics of the
polynomial $f$.

Together with the natural $\C$-linear structure, the relative
cohomology groups $\boldsymbol H^k_f$ possess the structure of a
module over the ring $\C[f]=f^*\C[x_1,\dots,x_n]$. This follows
from the identity
\begin{equation}\label{correctness}
  f\cdot(df\land\eta+d\theta)=df\land (f\eta-\theta)+d(f\theta).
\end{equation}
meaning that relatively exact forms are preserved by
multiplication by $f$.

As is well-known, the highest module $\boldsymbol H^n_f$, as well
as all $\boldsymbol H^k_f$ with $0<k<n-1$, is zero. Instead, we
consider another important module, called \emph{Brieskorn module}
(lattice) \cite{brieskorn,dimca-saito:brieskorn,douai}, defined as
the quotient
\begin{equation}\label{brieskorn}
  \boldsymbol B_f=\L^n/df\land d\L^{n-2},
\end{equation}
and the $\C[f]$-module $\boldsymbol P_f$, the quotient of
\emph{all} $(n-1)$-forms by the \emph{closed} $(n-1)$-forms,
\begin{equation}\label{petrov}
  \boldsymbol P_f=\L^{n-1}/(df\land\L^{n-2}+d\L^{n-2})
  \supseteq \boldsymbol H^{n-1}_f.
\end{equation}
The latter is an extension of $\boldsymbol H^{n-1}_f$: the
quotient $\boldsymbol P_f/H^{n-1}_f$ is naturally isomorphic to
the finite-dimensional $\C$-space $\L^n_f=\L^n/df\land \L^{n-1}$.
In several sources, $\boldsymbol P_f$ is referred to as the
\emph{Petrov module}. The exterior differential naturally projects
as a \emph{bijective} map $d\:\boldsymbol P_f\to \boldsymbol B_f$
which obviously \emph{is not} a $\C[f]$-module homomorphism.

Clearly, a relatively exact (closed) form is exact (resp., closed)
after being restricted on any nonsingular level set
$f^{-1}(t)\subset\C^n$, $t\in\C$, since $df$ vanishes on all such
sets.

The inverse inclusion is considerably more delicate. Gavrilov
studied the case $n=2$ and proved that for a 1-form with exact
restrictions on all level curves $f^{-1}(t)\subset\C^2$ to be
relatively exact, it is sufficient to require that the polynomial
$f$ has only isolated singularities and all level curves
$f^{-1}(t)$ be connected \cite{gavrilov:petrov-modules,
gavrilov:annalif}. This result generalizes the earlier theorem by
Ilyashenko \cite{ilyashenko:1969}. A multidimensional
generalization in the same spirit was obtained by I.~Pushkar$'$
\cite{arisha}. The affirmative answer depends on the topology of a
generic level set $f^{-1}(t)$ (its connectedness for $n=2$ or
vanishing of the Betti numbers $b_k$ for $k$ between $0$ and
$n-2$, see \cite{dimca-saito:general-fiber,bonnet-dimca}).

Both the isolatedness and connectedness assumptions can be derived
from a single assumption that the principal (quasi)homogeneous
part $\^f$ of the polynomial $f$ has an isolated critical point at
the origin: such polynomials are called
\emph{semiquasihomogeneous} \cite{odo-1}. For two variables with
equal weights it suffices to require that $\^f$ factors as a
product of pairwise different linear homogeneous terms.

\subsection{Computation of relative cohomology}
Besides the above question on the relationship between the
algebraically defined cohomology of the relative de~Rham complex
and analytically defined cohomology of (generic) fibers, the
natural problem of computing $\boldsymbol H^\bullet _f$ arises.

This problem was addressed in the papers
\cite{bonnet-dimca,dimca-saito:brieskorn,
dimca-saito:general-fiber,douai,gavrilov:annalif,gavrilov:petrov-modules}
mentioned above. Using analytic tools or theory of perverse
sheaves and $D$-modules, they prove that under certain
genericity-type assumptions on $f$, the highest relative
cohomology module $\boldsymbol H_f^{n-1}$ and the Petrov module
$\boldsymbol P_f$ are finitely generated over the ring $\C[f]$.
For semiquasihomogeneous polynomials one can describe explicitly
the collection of generators for $\boldsymbol B_f$, the polynomial
forms $\omega_1,\dots,\omega_l\in\L^{n-1}$ such that any other
form $\omega\in\L^{n-1}$ can be represented as
\begin{equation}\label{rce}
  \omega=\sum_{i=1}^l p_i\,\omega_i+df\land \eta+d\xi,
  \qquad p_i=p_i(f)\in\C[f],\ \eta,\xi\in\L^{n-2},
\end{equation}
with appropriate polynomial coefficients $p_i$ that are uniquely
defined.

The proofs of this and related results, obtained in either
analytic or algebraic way, are sufficiently involved. In
particular, it is very difficult if possible at all to get an
information on (i) how the decomposition \eqref{rce} depends on
parameters, in particular, if $f$ itself depends on parameters,
and (ii) how to place explicit \emph{quantitative bounds} on the
coefficients $p_i(f)$ in terms of the magnitude of coefficients of
the form $\omega$. For example, to extract such bounds from the
more transparent analytic proof by Gavrilov, one should place a
\emph{lower} bound on the determinant of the period matrix of the
forms $\omega_i$ over a system of vanishing cycles on the level
curves $f^{-1}(t)$. The mere nonvanishing of this determinant is a
delicate assertion whose proof in \cite{gavrilov:petrov-modules}
is incomplete (a simple elementary proof was supplied by Novikov
\cite{mit:irredundant}). The explicit computation of this
determinant for a specific choice of the generators $\omega_i$ was
achieved by A.~Glutsuk \cite{glutsuk}, but the answer is given by
a very cumbersome expression.

In the next section we suggest an elementary derivation of the
formula \eqref{rce} under the assumption that the polynomial $f$
is {semiquasihomogeneous}. This derivation:
\begin{enumerate}
  \item gives an independent elementary demonstration of the
  Gavrilov--Bonnet--Dimca theorem for the most important particular case
  of semiquasihomogeneous polynomials;
  \item proves that the polynomial coefficients $p_i$ and the forms
  $\eta$, $\theta$ from the decomposition \eqref{rce} depend polynomially
  on the coefficients of the non-principal part of $f$, provided that the
  principal quasihomogeneous part of $f$ remains fixed;
  \item yields the collection of the coefficients $(p_1,\dots,p_l)$ of
  \eqref{rce} as a result of application of a certain \emph{linear
  operator} to the form $\omega$. The norm of this operator can be
  explicitly bounded in terms of $f$ (and the chosen set of generators
  $\{\omega_i\}$) and the degree $\deg\omega$.
\end{enumerate}

\section{Bounded decomposition in the Brieskorn and Petrov
modules}

\subsection{Degrees, weights, norms}
In this section we first consider quasihomogeneous polynomials
from the ring $\C[x]=\C[x_1,\dots,x_n]$ with rational positive
weights $w_i=\deg x_i$ normalized by the condition
$w_1+\cdots+w_n=n$ to simplify the treatment of the most important
\emph{symmetric} case when $w_i=1$. The symbol $\deg f$ always
means the quasihomogeneous degree.

\begin{Rem}
Later on we will introduce additional variables
$\l=(\l_1,\dots,\l_m)$ considered as \emph{parameters}, assign
them appropriate weights and work in the extended ring
$\C[x,\l]=\C[x_1,\dots,x_n,\l_1,\dots,\l_m]$. Even in the
symmetric case the weights of the parameters will in general be
different from $1$.
\end{Rem}

The \emph{Euler field} associated with the weights $w_1,\dots,w_n$
is the derivation $X=\sum w_i\,x_i\partial/\partial x_i$ of
$\C[x]$. By construction, $Xf=r f$, $r=\deg f\in\mathbb Q$, for
any quasihomogeneous polynomial $f$ (the Euler identity).

We put $\deg dx_i=\deg x_i=w_i$. This extends the quasihomogeneous
grading on all $k$-forms: in the symmetric case, the degree of a
polynomial $k$-form will be $k$ plus the maximal degree of its
coefficients. Obviously, $\deg\omega=\deg d\omega$ for any form,
provided that $d\omega\ne0$. The Lie derivative $X\omega$ of a
quasihomogeneous form $\omega$ of degree $r$ by the Euler identity
is $r \omega$. Note that $\deg\omega>0$ for all $k$-forms with
$k\ge 1$.

The \emph{norm} of a polynomial in one or several variables is
defined as the sum of absolute values of its (real or complex)
coefficients. This norm is multiplicative. The norm of a $k$-form
by definition is the sum of the norms of its polynomial
coefficients; it satisfies the inequality
$\|\omega\land\eta\|\le\|\omega\|\cdot\|\eta\|$ for any two forms
$\omega,\eta$.

The exterior derivative operator is bounded in the sense of this
norm if the degree is restricted: $\|d\omega\|\le (\max_i
w_i)\deg\omega\cdot \|\omega\|$. In particular, in the symmetric
case $\|d\omega\|\le r\,\|\omega\|$, $r=\deg\omega$. Conversely, a
primitive of an $n$-form $\mu$ can be always chosen bounded by the
same norm $\|\mu\|$.

Unless explicitly stated differently, a monomial (monomial form,
etc) has always the unit coefficient.

\subsection{Parameters}
We will systematically treat the case when all objects (forms,
functions etc.) depend polynomially on finitely many additional
parameters  $\l=(\l_1,\dots,\l_m)$. We will denote by $\L^k[\l]$,
$k=0,\dots,n$, the collection of $k$-forms whose coefficients
polynomially depend on $\l$. For instance, the notation
$\eta\in\L^{n-1}[\l]$ means that $\eta=\sum_{i=1}^n
a_i(x,\l)\,dx_1\land\cdots\land\widehat{dx_i}\land\cdots\land
dx_n$ with polynomial coefficients $a_i\in\C[x,\l]$.

In such case the norm of forms, functions etc.\ will be always
considered relative to the ring $\C[x,\l]$, that is, as the sum
$\sum_i\|a_i\|$ of absolute values of coefficients $a_i$ of the
complete expansion in $x,\l$. If the parameters $\l_s$ are
assigned weights, we take them into account when defining the
degree of the form. To stress the fact that the norm is computed
relative to the ring $\C[x,\l]$ and not to $\C[x]$ (i.e., that the
situation is parametric), we will sometimes denote the norm by
$\|\cdot\|_\l$. For an instance, $\|2\l_1x_1\|=2|\l_1|\ne
2=\|2\l_1x_1\|_\l$.

\subsection{Division by a quasihomogeneous differential $df$.
The division modulus} If $f\in\C[x_1,\dots,x_n]$ is a
quasihomogeneous polynomial having an \emph{isolated} singularity
at the origin, then the multiplicity $l$ of this singularity can
be easily found by B\'ezout theorem, since no roots of the system
of algebraic equations $\partial f/\partial x_i=0$, $i=1,\dots,n$,
can escape to infinity. In the symmetric case $l=(\deg f-1)^n$.
Choose any monomial basis $\f_1,\dots,\f_l$ of the local algebra
$\C[[x_1,\dots,x_n]]/\<\partial f\>$, $\<\partial f\>=\big<\pd
f{x_1},\dots,\pd f{x_n}\big>$. Then the monomial $n$-forms
$\mu_i=\f_i\,dx_1\land \cdots\land dx_n$ form a basis of
$\L^n_f=\L^n/df\land \L^{n-1}$ over $\C$: any $n$-form $\mu$ can
be divided out as
\begin{equation}\label{divbydf}
  \mu=\sum_{i=1}^l c_i\mu_i+df\land\eta,
  \qquad c_i\in\C,\quad \eta\in\L^{n-1},
\end{equation}
with appropriate constants $c_1,\dots,c_l\in\C$ \(coefficients of
the ``remainder'' $\sum c_i\mu_i$\) and a polynomial form
$\eta\in\L^{n-1}$ \textup{(}the ``incomplete ratio''\textup{)}.
Moreover, if $\mu$ is quasihomogeneous, then the decomposition
\eqref{divbydf} contains only terms with $\deg\mu_i=\deg\mu$ and
$\deg\eta=\deg\mu-\deg f$. This immediately follows from
quasihomogeneity and the uniqueness of the coefficients $c_i$.
Form this observation we also conclude that all monomial forms of
degree $<\deg f$ must be among $\mu_i$, and, moreover, any
monomial form of degree greater than $\max_i\deg\mu_i$, is
divisible without remainder by $df$.

The choice of the monomial forms $\mu_i$ spanning the quotient, is
not unique, though the distribution of their degrees is. Denote by
$\rho=\rho(f)$ the maximal difference
\begin{equation}\label{rho}
  \rho(f)=\max_i\deg\mu_i-\min_i\deg\mu_i=\max_i\deg\f_i-\min_i\deg\f_i.
\end{equation}
The following results are well-known.

\begin{Prop}\label{prop:odo}\indent
1. In the symmetric case $\rho(f)<l=(r-1)^n$ \cite[\paragraph
5.5]{odo-1}.

2. In the bivariate case $n=2$ the inequality $\rho(f)<r=\deg f$
holds if and only if $f$ is a ``simple singularity'' of one of the
following types,
\begin{equation*}
  \begin{aligned}
  A_k:\quad& f=x_1^{k+1}+x_2^2,\qquad k\ge 2,
  \\
  D_k:\quad&f=x_1^2x_2+x_2^{k-1},\qquad k\ge 4,
  \\
  E_6:\quad&f=x_1^3+x_2^4,
  \\
  E_7:\quad &f=x_1^3+x_1x_2^3,
  \\
  E_8:\quad&f=x_1^3+x_2^5,
  \end{aligned}
\end{equation*}
see e.g., \cite[\paragraph 13, Theorem~2]{odo-1}.\qed
\end{Prop}

From these observations it can be immediately seen that the
division with remainder \eqref{divbydf} is a bounded linear
operation in the space of all $n$-forms of restricted degrees.

\begin{Lem}\label{lem:divbyhomogen}
Assume that $f\in\L^0$ is a quasihomogeneous polynomial having an
isolated critical point of multiplicity $l$ at the origin, and the
monomial $n$-forms $\mu_1,\dots,\mu_l\in\L^n$, form the basis of
$\L^n_f$.

Then there exists a finite constant $M<+\infty$ depending only on
$f$ and the choice of the basis $\{\mu_i\}$, such that any
$n$-form $\mu\in\L^n$ can be divided with remainder by $df$ as in
\eqref{divbydf} subject to the following constraints,
\begin{equation}\label{M(f)}
  \deg\eta\le\deg\mu-\deg f,\qquad \|\eta\|+\sum|c_i|\le M\|\mu\|.
\end{equation}
If the form $\mu$ is quasihomogeneous, then $\deg\eta=\deg\mu-\deg
f$ and $c_i$ can be nonzero only if $\deg\mu_i=\deg\mu$.
\end{Lem}

The constant $M$ depends on the choice of the monomial basis
$\{\mu_i\}$. The optimal choice of such basis (out of finitely
many possibilities) results in the smallest value $M=M(f)$ that
depends only on $f$. We will always assume that the basis
$\{\mu_i\}$ is chosen optimal in this sense.

\begin{Def}
The minimal constant $M(f)$ corresponding to an optimal choice of
the monomial basis of the quotient $\L^n_f$, is called the
\emph{division modulus} of the quasihomogeneous polynomial
$f\in\L^0$.
\end{Def}

\begin{Cor}
Assume that $\mu\in\L^n[\l]$ depends polynomially on additional
parameters $\l$. Then $\mu$ can be divided with remainder by $df$
so that the remainder and the incomplete ratio depend polynomially
on $\l$ with the same division modulus,
\begin{gather*}
  c_i=c_i(\l)\in\C[\l],\quad i=1,\dots,n,\qquad
  \eta\in\L^{n-1}[\l],
  \\
  \|\eta\|+\sum\|c_i\|\le M(f)\,\|\mu\|,\qquad
  \|\cdot\|=\|\cdot\|_{\l}.
\end{gather*}
\end{Cor}

\begin{proof}[Proof of the corollary]
Every monomial from the expansion of $\mu$ in $x,\l$ can be
divided out separately by $df$ which is independent of $\l$.
\end{proof}

\begin{proof}[Proof of the Lemma]
Let $M$ be the best constant such that \eqref{M(f)} holds for all
monomial $n$-forms with $\deg\mu\le l$. It is finite since there
are only finitely many such forms. In particular, since any form
of degree $l$ is divisible by $df$ by Proposition~\ref{prop:odo},
the respective fraction $\eta$ will be of the norm at most
$M\|\mu\|$.

Writing an arbitrary monomial $n$-form of degree $> l$ as a
product of a monomial form of degree $l$ times a monic monomial
function $x^\alpha\in\C[x]$, $\alpha\in\mathbb Z^n_+$, we
construct the explicit division formulas (without remainders) for
all monomial forms of higher degrees. The division constant will
be given by the same number $M$, since multiplication by a monic
monomial preserves the norms of both $\|\mu\|$ and $\|\eta\|$.

All the other assertions of the Lemma are well-known \cite{odo-1}.
\end{proof}

\subsection{Computability of the division modulus}
Despite its general nature, the above proof is constructive, at
least in the low dimensional cases $n=1,2$, allowing for an
explicit computation of the division modulus in these cases.

The one-dimensional case is trivial: for the monomial $f(x)=x^{r}$
the division modulus $M(f)$ is equal to $r$ and it can be
obviously recalculated for any other principal homogeneous part.
The ``special case'' of a multivariate polynomial
$f(x)=x_1^{r}+\cdots +x_n^{r}$, see \cite{glu-il}, is reducible to
the one-dimensional situation. In this case $l=(r-1)^n$ monomial
forms $x^\alpha\,dx_1\land \cdots\land dx_n$ with $0\le
\alpha_i\le r-1$ form the basis, and the corresponding division
modulus is again equal to $r$. This example admits an obvious
generalization for quasihomogeneous ``special polynomials'' with
different weights.

For a bivariate truly homogeneous polynomial $f$ (i.e., in the
symmetric case, the  most important for applications), the
division modulus $M$ for all higher degree forms ($\deg\mu\ge2\deg
f$) can be explicitly computed as the norm of the inverse
Sylvester matrix for the partial derivatives $\pd{f}{x_1}$ and
$\pd f{x_2}$ \cite{redundant}. The ``\emph{quasimonic}''
polynomials, introduced in that paper, are defined by the
condition $M(f)=1$, which in many respects is a natural
normalizing condition for multivariate polynomials.

The choice of the basic forms even in the symmetric bivariate case
depends on $f$: while it is generically possible to choose them as
$x_1^{\alpha_1}x_2^{\alpha_2}\,dx_1\land dx_2$ with $0\le
\alpha_{1,2}\le r-1$, for a badly chosen $f$ some of these forms
of degree greater than $r=\deg f$ can become linear dependent
modulo $df$, requiring a different choice. In order to avoid
making this choice, one may allow a \emph{redundant} (i.e., linear
dependent) collection of generating forms $\mu_i$. Choosing all
monomial forms of degree $\le 2r$ makes the corresponding division
for low degree forms trivial, so that the  division modulus $M(f)$
is determined only by division of forms of higher degree. Details
and accurate estimates in the bivariate symmetric case can be
found in \cite{redundant}.

To describe the division modulus $M(f)$ in the case of $n\ge 3$
variables is a considerably more difficult problem, though it
still can be reduced to analysis of finitely many monomial
divisions. One can (at least, theoretically) express $M(f)$ via
lower bounds for minors of certain explicitly formed matrices.

\begin{Rem}\label{rem:norm-vs-M}
It is worth mentioning that the division modulus $M(f)$ is not
directly related to the norm $\|f\|$, even in the symmetric
bivariate case. If $\deg\mu\ge l$ and $\mu=df\land\eta$, then
$\|\mu\|\le\|df\|\,\|\eta\|$. On the other hand, $\|\mu\|\ge
M^{-1}\|\eta\|$ by the definition of $M(f)$. Therefore
$M(f)\ge\|df\|^{-1}\ge \|f\|^{-1}$, that is, the division modulus
for a polynomial $f$ with the small norm must be large. The
inverse is not true: a polynomial with a small division modulus
can have a very large norm. Simple examples can be constructed in
the form $f(x)=c\prod_i(x_1-\l_ix_2)$ with sufficiently close
values of the parameters $\l_i\in[0,1]$ and a suitably chosen
normalizing constant $c\in\C$.
\end{Rem}

\subsection{Division by $f$}
We begin by establishing an analog of the Euler identity in the
Brieskorn module. It plays the central role for explicitly
constructing the decomposition \eqref{rce}.

\begin{Lem}\label{lem:Euler}
Assume that $f\in\L^0$ is a quasihomogeneous polynomial of degree
$r$. Then any polynomial $n$-form divisible by $df$ in $\L^n$, can
itself  be divided by $f$ in the Brieskorn module $\boldsymbol
B_f$. It also admits a polynomial primitive divisible by $f$.

In other words, for any form $\eta\in\L^{n-1}$ there exist four
forms $\mu\in\L^n$, $\omega\in\L^{n-1}$ and $\xi,\xi'\in\L^{n-2}$
such that
\begin{align}
  df\land \eta&=f\mu+df\land d\xi\label{div-f1}
  \\
  &=d(f\omega)+df\land d\xi'.\label{div-f2}
\end{align}

The degrees of all forms $\mu,\omega,\xi,\xi'$ are all equal to
$\deg\eta$ in case the latter is quasihomogeneous.

The division operation is always well-posed in the sense that the
decomposition \eqref{div-f2} can be always chosen to meet the
inequality
\begin{equation}\label{Euler-norm}
  \|\omega\|+\|\xi'\|\le (n+3)\deg\eta\cdot\|\eta\|
\end{equation}
\(a similar inequality can be proved also for the first
decomposition \eqref{div-f1}\).
\end{Lem}

\begin{proof}
Note that for any $n$-form $\mu\in\L^n$ and any vector field $X$
on $\C^n$,
\begin{equation*}
  (Xf)\mu=(\ix df)\,\mu=df\land\ix\mu,
\end{equation*}
where $\ix$ is the inner antiderivative, since $df\land\mu=0$. We
will need this formula for the case when $X$ is the Euler vector
field.

To prove the first divisibility assertion \eqref{div-f1}, we have
to show that the identity
\begin{equation}\label{1}
  df\land  \eta=f\mu+df\land d\xi
\end{equation}
can be always resolved as a linear equation with respect to $\mu$
and $\xi$ for any choice of $\eta$. Using the Euler identity for
functions and the above remark, we represent $f\mu$ as a form
divisible by $df$,
\begin{equation}\label{by-euler}
  f\mu=r^{-1}(Xf)\,\mu=r^{-1}(\ix df)\mu=df\land r^{-1}\,\ix \mu.
\end{equation}
The equation \eqref{1} will obviously be satisfied if
\begin{equation*}
  \eta=r^{-1}\,\ix\mu+d\xi,
\end{equation*}
that is, when $\eta$ is cohomologous to $\ix\mu$. This last
condition is equivalent to the equality between the exterior
derivatives
\begin{equation*}
  d\eta=r^{-1}\,d\ix\mu=r^{-1}X\mu,
\end{equation*}
since by the homotopy formula, $d\ix \mu=X\mu-\ix d\mu=X\mu$. Thus
resolving the equation \eqref{1} is reduced to inverting the Lie
derivative $X$ on the linear space of $n$-forms.

We claim that the linear map $\mu\mapsto X\mu$ of $\L^n$ to
itself, is surjective (and obviously degree-preserving),
guaranteeing thus solvability of the last equation for any choice
of $\eta$. Indeed, any monomial $n$-form
$\mu_\alpha=x^\alpha\,dx_1\land \cdots\land dx_n$ is an
eigenvector of $X$ with the strictly positive eigenvalue
$\deg\mu_\alpha\ge n$ (recall that the weights $w_i$ are
normalized so that the volume form $dx_1\land\cdots\land dx_n$ is
of degree $n$). Thus $X$ is surjective on $\L^n$ (actually,
bijective) and one can choose $\mu=rX^{-1}(d\eta)$. The norm of
the inverse operator $X^{-1}$ does not exceed $(r/n)\deg\eta$ in
the symmetric case. The proof of \eqref{div-f1} is complete.

To prove the second assertion \eqref{div-f2}, we transform it
using \eqref{by-euler} as follows,
\begin{equation*}
  df\land \eta=f\,d\omega+df\land(\omega+d\xi')=r^{-1}df\land\ix
  d\omega+df\land(\omega+d\xi'),
\end{equation*}
which will be obviously satisfied if
\begin{equation}\label{for-xi}
  \eta=r^{-1}\,\ix d\omega+\omega+d\xi'.
\end{equation}
Taking the exterior derivative as before, we reduce this equation to the
form
\begin{equation*}
  d\eta=r^{-1}\,d\ix d\omega+d\omega=r^{-1}\,X\mu+\mu,\qquad \mu=d\omega.
\end{equation*}
Solvability of this equation with respect to $\mu$ (and hence to
$\omega$) for any left hand side $d\eta$ follows from
invertibility of the differential operator $r^{-1}X+\id$ on the
linear space of polynomial $n$-forms ($\id$ stands for the
identity operator). Exactly as in the previous situation, all
monomial $n$-forms are eigenvectors for $(r^{-1}X+\id)|_{\L^n}$
with the positive eigenvalues, all greater or equal to
$r^{-1}n+1$, hence $r^{-1}X+\id$ is invertible on $\L^n$ and
$\omega$ can be chosen as a primitive of
$(r^{-1}X+\id)^{-1}d\eta$.

To prove the inequality between the norms, notice that
$\mu=d\omega$ satisfies the inequality $\|\mu\|\le
\|d\eta\|\le\deg\eta\,\|\eta\|$. A primitive $\omega$ can be
always take of the norm $\|\omega\|\le\|d\omega\|$. Together this
yields $\|\omega\|\le\deg\eta\,\|\eta\|$.

The norm $\|\xi'\|$ can be found from \eqref{for-xi}. Clearly,
$\|\ix\mu\|\le n\|\mu\|$ because of the choice of the weights
$\deg x_i$ which satisfy the condition $\sum w_i=n$.  Substituting
this inequality into \eqref{for-xi}, we obtain
\begin{equation*}
  \|\xi'\|\le \|d\xi'\|\le \|\eta\|+n\|d\omega\|+\|\omega\|\le
(n+2)\deg\eta\,\|\eta\|,
\end{equation*}
since $\deg\omega=\deg\eta\ge1$.
\end{proof}

\subsection{Generating Petrov and Brieskorn modules: the
algorithm} Division by the gradient ideal together with the Euler
identity as formulated in Lemma~\ref{lem:Euler}, allows for a
constructive proof of the representation \eqref{rce} for an
arbitrary semiquasihomogeneous polynomial $F$.

Let $F=f+h\in\C[x_1,\dots,x_n]$ be a semiquasihomogeneous
polynomial with the principal quasihomogeneous part $f$ and the
lower-degree part $h$. Denote as before by
$\mu_1,\dots,\mu_l\in\L^n$ the forms spanning
$\L^n_f=\L^n/df\land\L^{n-1}$ (note that the quotient is computed
using only the principal part $f$). We claim that:
\begin{enumerate}
    \item any $n$-form $\mu\in\L^n$ can be represented as
    \begin{equation}\label{bri}
    \mu=\sum_{i=1}^l q_i\,\mu_i+dF\land d\zeta,\qquad
    q_i\in\C[F],\ \zeta\in\L^{n-2},
    \end{equation}
    \item any $(n-1)$-form $\omega\in\L^{n-1}$ can be represented
    as
    \begin{equation}\label{pet}
    \omega=\sum_{i=1}^l p_i\,\omega_i+dF\land\xi+d\xi',
    \qquad p_i\in\C[F],\ \xi,\xi'\in\L^{n-2}.
    \end{equation}
\end{enumerate}

The construction of the decomposition \eqref{bri} begins by
division of $\mu$ by $df$ as explained in
Lemma~\ref{lem:divbyhomogen}:
\begin{equation*}
  \mu=\sum c_i\mu_i+df\land\eta,\qquad c_i\in\C,\ \eta\in\L^{n-1}.
\end{equation*}
If $\deg\mu<r=\deg f=\deg F$, then the incomplete ratio is in fact
absent, $\eta=0$, and we arrive to a particular case of
\eqref{bri} with $q_i=c_i$ of degree $0$ (constants).

If $\deg\mu$ is higher than $r$, we transform the term $df\land
\eta$ using Lemma~\ref{lem:Euler} and then substitute $f=F-h$:
\begin{equation*}
\begin{aligned}
  \mu-\sum c_i\mu_i&=f\mu'+df\land d\zeta
  \\
  &=F\,\mu'+dF\land d\zeta-\mu'',
  \qquad \mu''=h\mu'+dh\land d\zeta.
\end{aligned}
\end{equation*}
Obviously, both $\mu'$ and $\mu''$ are of degree \emph{strictly
inferior} to $\deg \mu$, which allows to continue the process
inductively. Assuming that the representations \eqref{bri} are
known for both $\mu'$ and $\mu''$, we substitute them into the
last identity and after collecting terms arrive to a
representation for $\mu$. In the symmetric case the inductive
process cannot take more than $\deg\mu-r$ steps. It is a direct
analog of the process of division of univariate polynomials, see
also \cite{redundant}.

To construct \eqref{pet}, we divide $d\omega$ by $df$. If
$\deg\omega<r$, then the incomplete ratio is absent and we obtain
a special kind of \eqref{pet} exactly as before.

Otherwise in the division with remainder
\begin{equation*}
  d\omega=\sum_{i=1}^l c_i\,d\omega_i+df\land\eta,\qquad
  c_i\in\C[\l],\ \eta\in\L^{n-1}[\l],
\end{equation*}
substitute $df\land \eta=d(f\omega')+df\land d\xi$ and pass to the
primitives. We obtain
\begin{equation}\label{step-alg}
  \begin{aligned}
  \omega-\sum c_i\omega_i&=f\omega'+df\land \xi+d\xi'
  \\
  &=F\omega'+dF\land\xi+d\xi'-\omega'',\qquad
  \omega''=h\omega'+dh\land \xi.
  \end{aligned}
\end{equation}
For the same reasons as before, the degrees of $\omega',\omega''$
are strictly smaller than $\deg\omega$, hence the process can be
continued inductively.

\begin{Rem}
In a somewhat surprising way, it turned out impossible to
transform directly the decomposition \eqref{bri} for the form
$d\omega\in\L^n$ into \eqref{pet} for $\omega$.
\end{Rem}

\subsection{Effective decomposition in the Petrov module}
The construction above is so transparent that any qualitative as
well as quantitative assertion concerning these expansions, can be
immediately verified.

We will show that
\begin{enumerate}
 \item all terms of the decomposition \eqref{pet} depend
 polynomially on the lower order terms of $F$, assuming that the
 principal part if fixed, and
 \item the well-posedness of the construction is determined solely
 by the division modulus $M(f)$ of the principal homogeneous part.
\end{enumerate}

In order to formulate the result, consider a \emph{general
semiquasihomogeneous polynomial} with the prescribed principal
quasihomogeneous part,
\begin{equation}\label{F}
  F(x,\l)=f(x)+h(x,\l),\qquad
  h(x,\l)=\sum_{\deg f_s<\deg f}\l_s f_s(x),
\end{equation}
where $f_1,\dots,f_m\in\C[x_1,\dots,x_n]$ are all (monic)
monomials of degree strictly inferior to $r=\deg f$, arbitrarily
ordered. We treat the coefficients $\l_1,\dots,\l_m$ as the
parameters of the problem, assigning to them the weights so that
\begin{equation*}
  \deg\l_s+\deg f_s=\deg f=r\qquad\text{for all }s.
\end{equation*}
This choice makes the entire polynomial $F$ quasihomogeneous of
the same degree $r$ in the ring
$\C[x,\l]=\C[x_1,\dots,x_n,\l_1,\dots,\l_m]$. Instead of the ring
$\C[F]$, the coefficients $p_i$ of the decomposition \eqref{pet}
will belong to the ring $\C[F,\l]$ and their quasihomogeneity will
be understood in the sense that the formal variable $F$ is
assigned the weight $\deg F=r$.

\begin{Thm}\label{thm:gavrilov}
If the quasihomogeneous polynomial $f\in\C[x]$ has an isolated
critical point at the origin and $F\in\C[x,\l]$ is a general
semiquasihomogeneous polynomial \eqref{F}, then any polynomial
quasihomogeneous $(n-1)$-form $\omega\in\L^{n-1}[\l]$ of degree
$k$ can be represented as
\begin{equation}\label{rce-again}
  \omega=\sum_{i=1}^l p_i\,\omega_i+dF\land \xi+d\xi'.
\end{equation}

The coefficients $p_i\in\C[F,\l]$ and the $(n-2)$-forms
$\xi,\xi'\in\L^{n-2}[\l]$ are all polynomial and quasihomogeneous
jointly in $F,\l$ \(resp., in $x,\l$\) of the degrees
$k-\deg\omega_i$, $k-r$ and $k$ respectively.

The norm of the coefficients relative to the ring
$\C[F,\l_1,\dots,\l_m]$ is explicitly bounded in terms of $n,r,k$
and the division modulus $M(f)$. In particular, for the symmetric
case when $\deg x_1=\cdots=\deg x_n=1$,
\begin{equation}\label{all-bd}
  \sum_{i=1}^l\|p_i\|\le k!\,r^{k(n+3)}M^k\,\|\omega\|,\qquad
  k=\deg\omega,\ M=M(f),\ \|\cdot\|=\|\cdot\|_\l.
\end{equation}
\end{Thm}

\begin{Rem}
The fact that the form $\omega$ is quasihomogeneous, is not
important: any polynomial form is the sum of quasihomogeneous
parts, each of them being divisible separately.
\end{Rem}

\begin{Rem}
Even in the symmetric case, the degrees of the parameters are
different from $1$: $\deg\l_s=r-\deg f_s$ will take all natural
values from $1$ to $r$.
\end{Rem}

\begin{proof}[Proof of the Theorem]
The first assertion of the Theorem (on polynomiality and
quasihomogeneity) follows from direct inspection of the algorithm
described above, since all transformations on each inductive step
(exterior differentiation, division by $df$ which is independent
of $\l$, and the Euler identity in $\boldsymbol P_f$) respect the
quasihomogeneous grading.

The only assertion that has to be proved is that on the norms. In
order for a sequence of increasing with $k$ real constants $C_k>0$
to be upper bounds for the decomposition \eqref{rce-again},
\begin{equation*}
  \sum_{i=1}^l \|p_i\|\le C_k\,\|\omega\|,\qquad
  \text{for all $\omega$ with}
  \deg\omega\le k,
\end{equation*}
they should satisfy a certain recurrent inequality which we will
instantly derive from the suggested algorithm.

Denote by $p_i\in\C[F,\l]$ (resp., by $p_i'$ and $p_i''$) the
polynomial coefficients of the decomposition of the forms $\omega$
(resp., $\omega'$ and $\omega''$) from the identity
\eqref{step-alg}: since the degrees of both $\omega',\omega''$ are
less than $k$ and the sequence $C_k$ is increasing, we have
\begin{equation*}
  \sum_i \|p_i'\|\le C_{k-1}\|\omega'\|,\qquad
  \sum_i \|p_i''\|\le C_{k-1}\|\omega''\|.
\end{equation*}
Multiplication by $F$ corresponds to a shift of coefficients in
the decomposition of $\omega'$. Thus from \eqref{step-alg} follows
the inequality
\begin{equation*}
  \sum_i\|p_i\|\le\sum_i\|c_i\|+\sum_i\|p_i'\|+\sum_i\|p_i''\|\le
  \sum_i \|c_i\|+C_{k-1}(\|\omega'\|+\|\omega''\|).
\end{equation*}
By Lemma~\ref{lem:Euler}, $\|\omega'\|\le (n+3)k\,\|\eta\|$. The
norm of the inferior part $h$ is by definition equal to the number
of terms, that is, the number of monomials in $n$ variables of
degree $\le r-1$. Therefore $\|h\|\le r^n$ and $\|dh\|\le
r^{n+1}$. This implies an upper bound for $\|\omega''\|$:
\begin{equation*}
  \|\omega''\|\le \|h\|\,\|\omega'\|+\|dh\|\,\|\xi\|\le
  (\|h\|+\|dh\|)(\|\omega'\|+\|\xi\|)\le 2r^{n+1}(n+3)k\,\|\eta\|
\end{equation*}
by Lemma~\ref{lem:Euler}. Finally, $\|\eta\|+\sum\|c_i\|\le
M\,\|\omega\|$ by definition of the division modulus $M=M(f)$.
Assembling all these bounds together, we conclude that
\begin{equation*}
  \sum\|p_i\|\le M \|\omega\|+C_{k-1}\cdot 3r^{n+1}(n+3)k
  \,\|\omega\|.
\end{equation*}
Thus the increasing sequence $C_k\ge 1$ will form upper bounds for
the norms of the coefficients of decomposition for polynomial
forms of degree $\le k$, provided that
\begin{equation*}
  C_k\ge A k C_{k-1}, \qquad A\ge 4r^{n+1}(n+3)M\ge r^{n+3}M
\end{equation*}
(notice that $r\ge 2$), which can be immediately satisfied if we
put
\begin{equation*}
  C_k=k!\,r^{k(n+3)}M^k.
\end{equation*}
This proves the inequality for the norms.
\end{proof}

Note that the bound established in this Theorem, is polynomial in
$M=M(f)$ and (for a fixed $r$) factorial in $k=\deg\omega$, that
is, only slightly overtaking the exponential growth.

\subsection{Nonhomogeneous division}
By a completely similar procedure one can describe the result of
division by a \emph{nonhomogeneous} differential $dF$ as a
sequence of divisions by the principal homogeneous part $df$.

More precisely, if $\mu\in\L^n[\l]$ is a polynomial $n$-form
polynomially depending on the parameters $\l_1,\dots,\l_m$ and
$F=f+\sum \l_s f_s$ is as in \eqref{F}, then there exists a
representation
\begin{equation}\label{div-by-dF}
  \mu=\sum_{i=1}^l c_i(\l)\mu_i+dF\land \eta,\qquad
  c_1,\dots,c_n\in\C[\l],\quad \eta\in\L^{n-1}[\l],
\end{equation}
polynomially depending on parameters. If $\mu$ is
quasihomogeneous, then so are $c_i$ and $\eta$, with $\deg
c_i=\deg\mu-\deg\mu_i$ and $\deg\eta=\deg\mu-\deg F$. Moreover,
the ratio $(\|c_i\|_\l+\|\eta\|_\l)/\|\mu\|_\l$ is bounded in
terms of $\deg\mu$ and the division modulus $M(f)$ of $f$ only.

Indeed, dividing $\mu$ by $df$ yields
\begin{equation*}
  \mu=\sum c_i\mu_i+df\land\eta=\sum c_i\mu_i+dF\land\eta -\mu',\qquad
  \mu'=dh\land\eta,
\end{equation*}
where $h=F-f$, hence $\deg h<\deg f=\deg F$ and therefore
$\deg\mu'<\deg\mu$. This means that the process of division can be
continued inductively. Since $\|\mu'\|\le \|h\|\,\|\eta\|\le
\const_{r,n} M(f)\,\|\mu\|$, the norms of the remainder and the
incomplete ratio are bounded in terms of $M(f)$ and the degrees.
In the symmetric case the bound looks especially simple.

\begin{Prop}\label{prop:div-by-dF}
In the symmetric case of all weights equal to $1$, the division of
a form of degree $k=\deg\mu$ is bounded as follows,
\begin{equation*}
  \|\eta\|_\l+\sum_{i=1}^l\|c_i\|_\l\le M_k(F)\cdot \|\mu\|,
  \qquad M_k(F)= kr^{n(k-r)}\,(M(f))^k.
\end{equation*}
\end{Prop}

\begin{proof}
In this case $\|h\|\le r^n$, so that $\|\mu'\|\le Mr^n\|\mu\|$,
and finally $\|\eta\|+\sum\|c_i\|\le
M\|\mu\|(1+K+\cdots+K^{\deg\mu-r})$, where $K=Mr^n$. Thus the norm
of the non-homogeneous division operator obviously does not exceed
$M^k (kr^{n(k-r)})$. This expression is exponential in $k=\deg\mu$
and polynomial in $M=M(f)$.
\end{proof}

\section{Picard--Fuchs system for Abelian
integrals}\label{sec:picfuc}

Consider a quasihomogeneous polynomial $f\in\L^0$ of degree
$r=\deg f$ having an isolated singularity of multiplicity $l$ at
the origin. As before, let $\mu_1,\dots,\mu_l$ be generators of
$\L^n_f$ over $\C$ and $\omega_1,\dots,\omega_l$ their monomial
primitives. Consider the general semiquasihomogeneous polynomial
$F=f+\sum_1^m \l_s f_s\in\C[x,\l]$ as in \eqref{F} with the fixed
principal part $f$, whose coefficients $\l_1,\dots,\l_m$ are the
natural parameters. Consider in the parameter space $\C^m$ the
locus $\S$ such that for $\l\in\C^m\ssm\S$ the level set
$\{x\in\C^n\:F(x,\l)=0\}$ is a nonsingular algebraic hypersurface.
Denote by $\G=\G(\l)$, $\l\notin\S$, any continuous family of
$(n-1)$-cycles on the zero level. The Abelian integrals
\begin{equation}\label{ai}
  I_i(\l)=\int_{\G(\l)}\omega_i,\qquad i=1\dots,l
\end{equation}
are well defined multivalued analytic functions on $\C^m\ssm\S$.
In this section we will derive a Pfaffian system of linear
equations satisfied by these integrals.

We will always assume that the weights of the parameters $\l_s$
are chosen so that $F$ becomes a quasihomogeneous polynomial in
$x,\l$ of degree $r$: $\deg\l_s=r-\deg f_s$. The enumeration of
the monomials $f_s$ begins with the free term $f_1\equiv1$ of
degree $0$ so that the respective coefficient $\l_1$ is
necessarily of degree $r$. Recall that $\rho(f)$ is the maximal
difference \eqref{rho} between the degrees of the forms $\mu_i$.

\begin{Thm}\label{thm:ppf}
There exist $(l\times l)$-matrix polynomials
$C_0(\l),C_1(\l),\dots,C_m(\l)$,
\begin{equation}\label{C}
\begin{gathered}
   C_0(\l)=\l_1\cdot \id+C'(\l_2,\dots,\l_m),
   \\
   \deg C_0\le r+\rho(f),\quad \deg C_s\le\deg f_s+\rho(f),
   \qquad s=1,\dots,m
\end{gathered}
\end{equation}
\(the degrees are quasihomogeneous\), such that on $\C^m\ssm\S$
\begin{equation}\label{pd}
  \pd{}{\l_s}\bigl(C_0(\l)I\bigr)=C_s(\l)I,\qquad s=1,\dots,m.
\end{equation}

The norms $\|C_s\|_\l$ are bounded by a power of the division
modulus $M(f)$.
\end{Thm}

In other words, the column vector function $I(\l)$ on the
complement to $\S$ satisfies the matrix Pfaffian equation
\begin{equation}\label{ppf}
  \mathbf d I=\Omega I,\qquad \Omega=C_0^{-1}\cdot\bigg(-\mathbf d
  C_0+\sum_{s=1}^m C_s\,d\l_s\bigg),
\end{equation}
with a rational matrix-valued 1-form $\Omega$ having the poles
only on the locus $\S'=\{\det C_0=0\}\subset\C^m$. Here $\mathbf
d$ is the exterior derivation with respect to the variables $\l_s$
only: for $c(\l)\in\C[\l]$, $\mathbf dc=\sum_s
\pd{c(\l)}{\l_s}\,d\l_s$.

The proof is constructive. The description of the matrix
polynomials $C_s(\l)$ is given below.

\subsection{Gelfand--Leray derivative with respect to parameters}
\begin{Lem}\label{lem:gel-ler}
If $\omega\in\L^{n-1}$ is a polynomial form with constant
\(independent of $\l$\) coefficients, and $\eta_s\in\L^{n-1}[\l]$
any form satisfying the identity
\begin{equation}\label{gl-par}
  f_s\,d\omega=-dF\land \eta_s,
\end{equation}
\(recall that $f_s=\pd{F}{\l_s}$\), then
\begin{equation*}
  \pd{}{\l_s}\int_{\G(\l)}\omega=\int_{\G(\l)}\eta_s.
\end{equation*}
\end{Lem}

\begin{proof}
To derive this formal identity, we express $\l_s=H(x)$ from the
equation $F(x,\l_s)=0$, assuming all other parameters fixed, and
apply the Gelfand--Leray formula to $H$: for \eqref{gl-par} to
hold, it would be sufficient if $\eta=\eta_s$ satisfies
\begin{equation*}
  d\omega=dH\land\eta.
\end{equation*}
It remains to observe that by the implicit function theorem and
the definition of the parameters,
\begin{equation*}
  dF+\pd{F}{\l_s}dH=0,\qquad \pd{F}{\l_s}=f_s.
\end{equation*}
Here and above $d$ stands for the exterior derivative with respect
to the ``spatial'' variables $x_1,\dots,x_n$.
\end{proof}

The standard Gelfand--Leray derivative appears for the parameter
occurring before the constant term $f_1\equiv 1$ (modulo the
sign).

\subsection{Derivation of the system: beginning of the proof
of Theorem~\ref{thm:ppf}}
Divide each of the forms $F\mu_i\in\L^n[\l]$, $\mu_i=d\omega_i$,
by $dF$ with with the remainder coefficients and the incomplete
ratios polynomially depending on $\l$ as in
Proposition~\ref{prop:div-by-dF}:
\begin{equation}\label{F/dF}
  F\mu_i=dF\land \eta_i+\sum_{j=1}^l c_{ij}\,\mu_j,\qquad
  c_{ij}=c_{ij}(\l).
\end{equation}
Clearly, the quasihomogeneous degree $\deg c_{ij}$ in $\C[\l]$ is
equal to $r+\deg \mu_i-\deg m_j\le \rho(f)+r$ ($c_{ij}\equiv0$ if
the difference is negative).

Let $C_0=C_0(\l)$ be the $(l\times l)$-matrix polynomial with the
entries $c_{ij}(\l)$. Since $dF$ does not depend on $\l_1$ (the
free term of $F$), while the only term depending on $\l_1$ in
$F\mu_i$ is $\l_1\mu_i$, the dependence of $C_0$ on $\l_1$ can be
immediately described: the corresponding remainder coefficients
$c_{ij}(\l_1)$ for the division of $\l_1\mu_i$ by $dF$ form the
scalar matrix $\l_1\cdot \id$ (the incomplete ratio is absent).

Since $c_{ij}$ do not depend on $x$ (being ``constants depending
on the parameters''), the identity \eqref{F/dF} implies that
\begin{equation*}
  d\bigl(F\omega_i-\sum\nolimits_j
  c_{ij}\omega_j\bigr)=-dF\land(-\omega_i-\eta_i),\qquad i=1,\dots,l.
\end{equation*}
Let
\begin{equation*}
  \omega'_{i,s}=-f_s(\omega_i+\eta_i),\qquad i=1,\dots,l,\quad
  s=1,\dots,m.
\end{equation*}
All these forms are polynomial and polynomially depending on
parameters. Their degrees can be easily computed:
$\deg\eta_i=\deg\mu_i=\deg\omega_i$, $\deg\omega'_{i,s}=\deg
f_s+\deg\mu_i$.

By the parametric Gelfand--Leray formula
(Lemma~\ref{lem:gel-ler}), the partial derivatives of integrals of
the forms $F\omega_i-\sum_j c_{ij}\omega_j$ over the cycle
$\G(\l)\subset\{F=0\}\subset\C^n$ are equal to the integrals of
the forms $\omega'_{i,s}$. Since the terms $F\omega_i$ vanish on
$\G(\l)$ for all values of $\l$, we have
\begin{equation*}
  \pd{}{\l_s}\biggl(\sum_j c_{ij}(\l)\,I_j(\l)\biggr)
  =I'_{i,s}(\l), \qquad I'_{i,s}(\l)=\oint_{\G(\l)}\omega'_{i,s}.
\end{equation*}
The forms $\omega_i$ were chosen to generate the Petrov module
$\boldsymbol P_F$ over $\C[F,\l]$, so each of the Abelian
integrals $\oint\omega'_{i,s}$ can be expressed as a polynomial
combination,
\begin{equation*}
  I'_{i,s}=\sum_{j=1}^l p_{ij,s}\,I_j,
  \qquad p_{ij,s}\in\C[F,\l],
\end{equation*}
for all $i,s$. Denote by $C_s=C_s(\l)$ the polynomial $(l\times
l)$-matrix function formed by the free terms of the polynomials
$p_{ij,s}(\cdot,\l)$:
\begin{equation*}
  C_s(\l)=\left[\left.p_{ij,s}(F,\l)\right|_{F=0}\right]_{i,j=1}^l,\qquad
  s=1,\dots,m.
\end{equation*}
All other terms, being divisible by $F$, disappear after
integration over the cycle on the level surface $\{F=0\}$.
Collecting the terms, we conclude that the partial derivatives of
the column vector function $I(\l)=(I_1(\l),\dots,I_l(\l))$,
$I_i=\oint \omega_i$, we have
\begin{equation*}
  \pd{(C_0I)}{\l_s}=C_sI,\qquad s=1,\dots,m.
\end{equation*}

\subsection{Bounds for the norms: end of the proof of
Theorem~\ref{thm:ppf}} The construction described above, does not
yet imply the assertion on the norms of the matrix polynomials
$C_0,\dots,C_m$ for only one reason: multiplication by $F=f+h$,
$h=\sum \l_s f_s$, is \emph{not} a bounded operator. While
multiplication by $h$ increases the norm at most by
$\|h\|_\l=\const_{n,r}$ (not exceeding $(r-1)^n$ in the symmetric
case), the norm $\|f\|$ cannot be bounded in terms of $M(f)$, as
required in the Theorem (see Remark~\ref{rem:norm-vs-M}).

To correct this drawback, exactly as in \cite{redundant}, the
division line \eqref{F/dF} should be first prepared using
\eqref{by-euler} as follows,
\begin{equation}\label{df/f}
\begin{gathered}
  F\mu_i=(f+h)\mu_i=df\land \eta'_i+h\mu_i=dF\land \eta_i'+\mu_i',
  \\
  \eta_i'=r^{-1}\ix\mu_i,\qquad \mu_i'=h\mu_i-dh\land \eta_i',
\end{gathered}
\end{equation}
where (we again make all estimates for the symmetric case only),
\begin{equation*}
  \|\eta_i'\|\le (n/r)\|\mu_i\|,\qquad
  \|\mu_i'\|\le\|h\|(1+r)(n/r)\|\mu_i\|.
\end{equation*}
Then forms $\mu_i'$ should be divided by $dF$ with remainder:
since their norms are bounded by a constant depending only on
$n,r$ (the norms of the monomial forms $\mu_i$ are equal to $1$),
the results of such division will be bounded by suitable powers of
$M(f)$ by virtue of Proposition~\ref{prop:div-by-dF}.

Collecting the terms, we conclude that the coefficients
$c_{ij}\in\C[\l]$ of the corresponding remainders in \eqref{F/dF}
and the incomplete ratios $\eta_i\in\L^{n-1}[\l]$ will be bounded
by expressions polynomial in $M(f)$.

The rest of the derivation remains unchanged and the estimates
completely straightforward: the polynomial bounds for $\eta_i$
imply those of the polynomial coefficients $p_{ij,s}\in\C[F,\l]$
by Theorem~\ref{thm:gavrilov}. This proves the last assertion of
Theorem~\ref{thm:ppf}. \qed

\section{Observations. Discussion}

The algorithm of derivation of the Picard--Fuchs system in the
Pfaffian form is so transparent that many things become obvious.

\subsection{Bounds} Though the matrix polynomials $C_s(\l)$ are not
quasihomogeneous (their entries have different degrees), the
determinant $\det C_0(\l)$ is a quasihomogeneous polynomial from
$\C[\l]$. Its degree can be immediately computed as $lr$ from the
explicit representation \eqref{C}. This same representation proves
that this determinant, equal to $\l_1^n+\text{polynomial
in}(\l_2,\dots,\l_m)$, does not vanish identically, so that the
system \eqref{ppf} is indeed meromorphic.

Moreover, the norm of the inverse matrix $C_0^{-1}$ can be
explicitly majorized in terms of the distance to the critical
locus. One possibility to do this is to consider the sections
$\l_1=1$ and apply the Cartan inequality as in \cite{redundant},
using the quasihomogeneity.

\subsection{Spectrum} The spectrum of $C_0(\l)$ can be also easily computed: it
consists of all $l$ critical values of the polynomial $F(x,\l)$,
at least when $F(\cdot,\l)$ is a Morse polynomial. To see this, it
is sufficient to evaluate both parts of \eqref{F/dF} at any of $l$
critical points $a_1,\dots,a_l\in\C^n$. The column vectors
$v_i=(\f_1(a_i),\dots,\f_l(a_i))^{T}$, $i=1,\dots,l$, are the
corresponding eigenvectors (recall that
$\mu_i=\f_i\,dx_1\land\cdots\land dx_n$).

\subsection{Hypergeometric form}
Restricting the Pfaffian system \eqref{ppf} on the
one-dimensional complex lines $\l_s=\const$, $s=2,\dots,m$,
parameterized by the value of $t=\l_1$, one obtains a
parameterized family of Picard--Fuchs systems of ordinary
differential equations. In this case only the matrix $C_1$ is
relevant.

By Theorem~\ref{thm:ppf}, it is quasihomogeneous of degree
$\le\rho(f)$ jointly in the variables $\l_1,\dots,\l_m$. If
$\rho(f)<r=\deg\l_1$, then $C_1$ cannot depend on $\l_1$ and hence
the Picard--Fuchs system in this case will have the hypergeometric
form \eqref{hyperg}. By Proposition~\ref{prop:odo}, this happens
only when $f$ is a simple quasihomogeneous polynomial of one of
the types listed there. For hyperelliptic polynomials (the
singularity of the type $A_k$) this was well-known, see
\cite{redundant}. In turn, the hypergeometric form implies that
all singular points of the Picard--Fuchs system are
\emph{Fuchsian}.

\subsection{Logarithmic poles}
For the full Pfaffian system \eqref{ppf} the polar locus,
occurring where $\det C_0(\l)$  vanishes, is of multiplicity $1$
(it is sufficient to produce just one value of the parameters $\l$
such that $F(\cdot,\l)$ has simple critical points). Yet it is not
the characteristic property.

A rational 1-form $\omega$ analytic outside a hypersurface
$\S'=\{g=0\}\subset\C^m$, $g$ being a polynomial without multiple
factors, is said to have a \emph{logarithmic singularity} on this
hypersurface, if both $g\omega$ and $dg\land\omega$ extend as
polynomial forms across $\S'$ on $\C^m$.

This is only one of several close but non-equivalent definitions,
probably the strongest possible. It ensures that the restriction
of $\omega$ on any holomorphic curve $\gamma$ cutting $\S'$ at a
point $a$, has a Fuchsian singularity with the residue independent
on the choice of $\gamma$, depending only on the point $a$.

The basic question concerning the system \eqref{ppf} is whether
this system itself or a suitable gauge transformation of this
system with a rational matrix gauge function, are Fuchsian with
bounded residues. If the answer is positive, this would mean a
positive solution of the infinitesimal Hilbert problem.

Using symbolic computation for implementing the algorithm, we
discovered that in the hyperelliptic case (singularity of the type
$A_k$) the Picard--Fuchs system \eqref{ppf} indeed has only
logarithmic poles until the degree $k=6$ of the polynomial
$f=x_1^k+x_2^2$. This naturally suggests the following conjecture.

\begin{Conj}
All singularities of the Picard--Fuchs system \eqref{ppf} are only
logarithmic poles on $\S'=\{\det C_0=0\}$.
\end{Conj}

It would be interesting to verify this conjecture for other simple
singularities listed in Proposition~\ref{prop:odo}, perhaps first
by symbolic computation.

The next step could be to study the behavior of residue of
\eqref{ppf}, the matrix function defined on the regular part of
$\S'$, checking whether it is bounded near singular points of the
discriminant.

\subsection{Singular perturbations}
The polynomial dependence of the matrices $C_s$ on the lower
degree coefficients of the polynomial $F=f+\cdots$ fails for the
coefficients of the principal part. Though apparently rational,
this dependence certainly must exhibit singularities when $f$
degenerates into a quasihomogeneous form with non-isolated
singularities. The Picard--Fuchs system in such cases may have
singular points corresponding to \emph{atypical values} of $F$.
Their appearance must somehow be related to the fact that the
division modulus explodes when such degeneracy occurs, thus
creating a \emph{singularly perturbed} system of linear
differential equations. These phenomena seem to be worth of
detailed study.


\def\BbbR{$\mathbf R$}\def\BbbC{$\mathbf
  C$}\providecommand\cprime{$'$}\providecommand\mhy{--}\font\cyr=wncyr9
\providecommand{\bysame}{\leavevmode\hbox
to3em{\hrulefill}\thinspace}
\providecommand{\MR}{\relax\ifhmode\unskip\space\fi MR }
\providecommand{\MRhref}[2]{%
  \href{http://www.ams.org/mathscinet-getitem?mr=#1}{#2}
} \providecommand{\href}[2]{#2}

\end{document}